\documentclass[11pt ,a4paper]{amsart}
\usepackage[T1,T2A]{fontenc}
\usepackage[utf8]{inputenc}
\usepackage[a4paper,margin=1in]{geometry}
\usepackage{fancyhdr}
\usepackage{graphicx}
\usepackage{pgf,tikz}
\usepackage{tikz-cd}
\usepackage{mathrsfs}
\usepackage{amsthm}
\usepackage{amssymb}
\usepackage{amsmath}
\usepackage{mathtools}
\usepackage{bm}
\usepackage{subfiles}
\usepackage{todonotes}
\usepackage{adjustbox}
\usepackage{hyperref}
\hypersetup{colorlinks=true, citecolor=black, linkcolor=black}
\usepackage[capitalize, nameinlink]{cleveref}
\usepackage{pdftexcmds}
\usepackage{bbm} 
\usepackage[bbgreekl]{mathbbol}
\usepackage[toc,page]{appendix}
\usepackage{enumitem}
\usepackage{xcolor}
\DeclareSymbolFontAlphabet{\mathbb}{AMSb}
\DeclareSymbolFontAlphabet{\mathbbl}{bbold}
\DeclareMathSymbol\bbDelta  \mathord{bbold}{"01}
\usepackage[backend=biber, style=alphabetic, maxnames=10, maxalphanames=3, minalphanames=3]{biblatex}
\usepackage[all]{xy}
\usetikzlibrary{arrows}
\theoremstyle{plain}

\makeatletter
\newcommand{\colim@}[2]{%
  \vtop{\m@th\ialign{##\cr
    \hfil$#1\operator@font colim$\hfil\cr
    \noalign{\nointerlineskip\kern1.5\ex@}#2\cr
    \noalign{\nointerlineskip\kern-\ex@}\cr}}%
}
\newcommand{\colim}{%
  \mathop{\mathpalette\colim@{\rightarrowfill@\textstyle}}\nmlimits@
}
\makeatother

\newtheorem{mthm}{Theorem}[section]
\newtheorem{mcor}[mthm]{Corollary}
\newtheorem{mlem}[mthm]{Lemma}
\newtheorem{mprop}[mthm]{Proposition}

\newtheorem*{mthm*}{Theorem}
\theoremstyle{definition}
\newtheorem{mrmk}[mthm]{Remark}
\theoremstyle{remark}

\theoremstyle{definition}
\newtheorem{mdef}[mthm]{Definition}

\theoremstyle{remark}

\title{The Universal Nilpotent Cone and A Canonical 1-Form}
\date{\today.}
\author{Marius Kjærsgaard}

\address{IHES, 35 route de Chartres, Bures-sur-Yvette 91440, France}

\email{kjaersgaard@ihes.fr}

\renewcommand{\varphi}{\phi}

\newcommand{\adjunction}[4]{\xymatrix@1{#1{\ } \ar@<0.3ex>[r]^{ {\scriptstyle #2}} & {\ } #3 \ar@<0.3ex>[l]^{ {\scriptstyle #4}}}}

\renewcommand{\epsilon}{\varepsilon}

\begin{document}
\begin{abstract}
    We construct a canonical Du Bois $1$-form on the nilpotent cone in the dual Lie algebra of any reductive group. Using this 1-form we give an intrinsic construction of Nadler and Yun's universal global nilpotent cone. In particular, we obtain a moduli interpretation of the universal nilpotent cone.
\end{abstract}

\maketitle

\setcounter{section}{0}
\section{Introduction}
Let $G$ be a split reductive algebraic group over a field $k$ of characteristic $0$, and $\mathrm{Bun}_G(\mathcal{C}_g/\mathcal{M}_g) $ be the moduli space of $G$-bundles on the universal curve $\mathcal{C}_g$ over the moduli space $\mathcal{M}_g$ of genus $g$ curves. 

\subsection{Nilpotent cones} The universal nilpotent cone $\mathrm{Nilp}_u\subset T^*\mathrm{Bun}_G(\mathcal{C}_g/\mathcal{M}_g)$ is a central object in the Betti geometric Langlands program for families of curves, see \cite{nadler2021automorphic}. However, an intrinsic characterization of universal nilpotence and a moduli interpretation have been missing. In this paper, we give a new construction of $\mathrm{Nilp}_u$ that provides both for smooth test schemes.

\subsubsection{The relative nilpotent cone} The relative nilpotent cone $\mathrm{Nilp}_\mathrm{rel}$ is a closed substack of the relative cotangent space $T^*_{\mathcal{M}_g}\mathrm{Bun}_G(\mathcal{C}_g/\mathcal{M}_g)$. A geometric point of $\mathrm{Nilp}_\mathrm{rel}$ can be characterized as a Higgs field $v\in H^0(X,\mathfrak{g}^*_{P_G}\otimes \Omega^1_X)$, for which there exist a $B$-reduction $P_B$ of $P_G$ such that $v$ is in the kernel of the map
\[H^0(X,\mathfrak{g}^*_{P_G}\otimes \Omega^1_X)\to H^0(X,\mathfrak{b}^*_{P_B}\otimes \Omega^1_X).\] 
This characterization does not immediately give a moduli interpretation of $\mathrm{Nilp}_\mathrm{rel}$, beyond the field valued points. Instead an intrinsic characterization of nilpotence, without reference to a $B$-reduction, can be given in terms of the relative Hitchin fibration. The relative Hitchin fibration defines a map
\[H_\mathrm{rel}:T^*_{\mathcal{M}_g}\mathrm{Bun}_G(\mathcal{C}_g/\mathcal{M}_g)\to \mathrm{Hitch}_{G,\mathrm{rel}},\]
from the cotangent space of $\mathrm{Bun}_G(\mathcal{C}_g/\mathcal{M}_g) $ relative to $\mathcal{M}_g$ to the relative Hitchin base. The relative nilpotent cone $\mathrm{Nilp}_\mathrm{rel}$ can then be defined as the preimage of the zero-section in the Hitchin base under the relative Hitchin fibration.

\subsubsection{The universal nilpotent cone}
Let $\pi:T^*\mathrm{Bun}_G(\mathcal{C}_g/\mathcal{M}_g)\to T^*_{\mathcal{M}_g}\mathrm{Bun}_G(\mathcal{C}_g/\mathcal{M}_g)$ denote the projection map from the absolute cotangent space to the cotangent space relative to $\mathcal{M}_g$. In \cite{nadler2021automorphic} and \cite{nadler2026global}, Nadler and Yun construct a Lagrangian lift of $\mathrm{Nilp}_\mathrm{rel}$ from the relative cotangent space to the absolute cotangent space. Their construction can be interpreted as follows: Fix $B\subset G$ a Borel and let $\pi^{-1}\mathrm{Nilp}_{\mathrm{rel},B}$ denote the preimage of the zero-section of $T^*_{\mathcal{M}_g}\mathrm{Bun}_B$ under the map 
\[\mathrm{Bun}_B\times_{\mathrm{Bun}_G} T^*\mathrm{Bun}_G\to T^*\mathrm{Bun}_B\to T_{\mathcal{M}_g}^*\mathrm{Bun}_B,\]
There is a map $H_{0,B}:\pi^{-1}\mathrm{Nilp}_{\mathrm{rel},B}\to T^*\mathcal{M}_g$ and the universal nilpotent cone is then defined as the image of $H_{0,B}^{-1}(\mathbf{0}_{T^*\mathcal{M}_g})$ under the projection $\pi^{-1}\mathrm{Nilp}_{\mathrm{rel},B}\to \pi^{-1}\mathrm{Nilp}_\mathrm{rel}$. 

This construction is analogous to the extrinsic description of the geometric points of $\mathrm{Nilp}_\mathrm{rel}$ in terms of $B$-reductions. In particular, this construction does not give a clear moduli interpretation of the universal nilpotent cone. Note that it is also not immediately clear from the construction that it defines a closed cone.

In this paper we give an intrinsic characterization of universal nilpotence in $T^*\mathrm{Bun}_G(\mathcal{C}_g/\mathcal{M}_g)$, similar to how the Hitchin fibration gives an intrinsic notion of nilpotence in the relative cotangent space not relying on $B$-reductions. This will give a new construction of the universal nilpotent cone, which in particular gives a clear moduli interpretation.

\subsubsection{A part of a universal Hitchin map}  The main construction of this paper is a canonical $\mathbb{G}_m$-equivariant morphism
\[H_0:\pi^{-1}\mathrm{Nilp}_{\mathrm{rel}}\to T^*\mathcal{M}_g,\]
which is uniquely characterized by the property that $H_{0,B}$ factors as the projection $\pi^{-1}\mathrm{Nilp}_{\mathrm{rel},B}\to \pi^{-1}\mathrm{Nilp}_\mathrm{rel}$ followed by $H_0$.
Thus we can define the universal nilpotent cone intrinsically as
\[\mathrm{Nilp}_\mathrm{u}:=H_0^{-1}(0_{T^*\mathcal{M}_g}).\]
Note that closedness is immediately clear from this construction. We prove the following properties of $H_0$, from which it follows that $\mathrm{Nilp}_u$ is Lagrangian.

\begin{mprop}
    The pullback of the tautological $1$-form on $T^*\mathcal{M}_g$ along $H_0$ agrees with the pullback of the tautological $1$-form on $T^*\mathrm{Bun}_G(\mathcal{C}_g/\mathcal{M}_g)$.
\end{mprop}
From this we deduce the following corollary.
\begin{mcor}
  The closed substack $H_0^{-1}(\mathbf{0}_{T^*\mathcal{M}_g})\subset T^*\mathrm{Bun}_G$ provides a conic Lagrangian lift of $\mathrm{Nilp}_\mathrm{rel}$.
\end{mcor}
The cone constructed this way agrees with the universal nilpotent cone defined by Nadler and Yun. The construction of $H_0$ is given in terms of moduli problems and thus we obtain a moduli description of $T$-valued points of $\mathrm{Nilp}_u$ for smooth test objects $T$.

\subsection{Du Bois forms and nilpotent cones in Lie algebras}
The map $H_0$ arises from a global Du Bois 1-form that we construct canonically on the nilpotent cone in the dual Lie algebra of any reductive group.

Let $\mathfrak{g}$ be the Lie algebra of $G$ and $\mathcal{N}\subset \mathfrak{g}^*$ denote the nilpotent cone in the dual of $\mathfrak{g}$. By a Du Bois 1-form we mean a section of $\Omega^{1,h}$, the sheafification of the Kähler 1-forms in the h-topology. It can be shown that Du Bois forms on a smooth variety identify with Kähler forms, see \cite{huber2014differentialformshtopology}. Let $\mathcal{B}$ denote the variety of Borel subgroups of $G$, we will prove that the tautological 1-form $\omega_\tau$ on $T^*\mathcal{B}$ descends to a global Du Bois 1-form $\omega_\mathcal{N}$ on $\mathcal{N}$ under the Springer resolution $\sigma:T^*\mathcal{B}\to \mathcal{N}$.

\begin{mprop}
    There exists a unique element $\omega_\mathcal{N}\in H^0(\mathcal{N},\Omega^{1,h}_\mathcal{N})$ such that $\sigma^*\omega_\mathcal{N}=\omega_\tau$ under the identification $H^0(T^*\mathcal{B},\Omega^{1,h}_{T^*\mathcal{B}})\cong H^0(T^*\mathcal{B},\Omega^1_{T^*\mathcal{B}})$
\end{mprop}

We emphasize that the sheafification in the $h$-topology is used twice in the proof, once in order to have descent along the Springer map. The second use of the h-sheafification takes care of the problem that $\omega_\tau$ does not vanish on the Springer fibers as a Kähler form, however, it does vanish on the Springer fibers as a Du Bois form.

    \subsection{Extension to $\operatorname{Spec } \mathbb{Z}$} We end the paper with a calculation for $\mathrm{PGL}_2$ over $\operatorname{Spec }\mathbb{Z}$, which indicates that our construction of $\mathrm{Nilp}_u$ has an analogue over a much more general base ring $k$ than fields of characteristic $0$. In the case $G=\mathrm{PGL}_2$, we prove an analogue of Proposition 1.3 over $\operatorname{Spec } \mathbb{Z}$. This allows us to give a construction of the universal nilpotent cone over the integers in this case. This direction does not seem directly relevant for the applications we have in mind, so we will not pursue it further than $\mathrm{PGL}_2$.

\subsection{Acknowledgements} I would like to thank Dustin Clausen for helpful discussions and Joakim Færgeman for a fruitful collaboration on related topics.

\section{Notation}
Throughout the text, unless otherwise stated, $k$ is a field of characteristic $0$ and $G$ is a split reductive group over $k$. All stacks will refer to sheaves of anima on the site $(\mathrm{Sm}/k)_{\mathrm{fppf}}$ of smooth affine schemes over $k$ with the fppf topology. By $G$-bundle we mean a left principal $G$-bundle.

\section{Du Bois Forms and Atiyah Bundles}

In this section we recall some basic properties of Du Bois forms and the Atiyah sequence.

\subsection{Du Bois forms} Du Bois forms are a modification of Kähler forms better adapted for certain problems involving singular varieties and recovers Kähler forms when applied to smooth varieties.

Let $\mathrm{Sch}/k$ denote the category of separated finite type $k$-schemes. Sending a scheme $X$ to its space of global Kähler 1-forms defines a functor
\[\Omega^1:(\mathrm{Sch}/k)^\mathrm{op}\to \mathrm{Vect}_k.\]
Let $\Omega^{1,h}:(\mathrm{Sch}/k)^\mathrm{op}\to \mathrm{Vect}_k$ denote the sheafification of $\Omega^1$ with respect to the h-topology on $\mathrm{Sch}/k$. 

\begin{mdef}
    For $X\in \mathrm{Sch}/k$, let $\Omega^{1,h}_{X/k}$ denote the sheaf on $X$ given by the restriction of $\Omega^{1,h}$ to the small Zariski site of $X$. 
\end{mdef}

By \cite[Proposition 4.2(i)]{huber2014differentialformshtopology}, the sheaf $\Omega^{1,h}_{X/k}$ is a torsion-free coherent sheaf of $\mathcal{O}_X$-modules on $X$. The natural map to the sheafification induces a map of coherent sheaves 
\[\Omega^1_{X/k}\to \Omega^{1,h}_{X/k}\] 
over $X$. The following proposition is \cite[Proposition 4.2(ii)]{huber2014differentialformshtopology}.

\begin{mprop}(\cite[Proposition 4.2(ii)]{huber2014differentialformshtopology}).
    If $X\in \mathrm{Sch}/k$ is smooth, then the natural map 
    \[\Omega^1_{X/k}\to \Omega^{1,h}_{X/k}\]
    is an equivalence.
\end{mprop}

Let $f:X\to Y$ be a morphism in $\mathrm{Sch}/k$. We will also need a relative version of Du Bois 1-forms.

\begin{mdef}
The sheaf of relative Du Bois differentials $\Omega^{1,h}_{X/Y}$ is the cokernel of the map
    \[f^*\Omega^{1,h}_{Y/k}\to \Omega^{1,h}_{X/k}.\]
    In particular there is a natural map $\Omega^1_{X/Y}\to \Omega^{1,h}_{X/Y}$.
\end{mdef}

If $X\to Y$ is a morphism between smooth varieties, then this agrees with the relative Kähler differentials. 

\subsection{Atiyah bundles} For a morphism $X/T$ and a $G$-bundle $p:P_G\to X$, the sequence $0\to p^*\Omega^1_{X/T}\to \Omega^1_{P_G/T}\to \Omega^1_{P_G/X}\to 0$ of coherent sheaves on $P_G$ descends along the $G$-action on $P_G$ to a sequence
\[0\to \Omega_{X/T}^1\to \mathrm{At}^*_{X/T}(P_G)\to \mathfrak{g}^*_{P_G}\to 0\]
on $X$. We refer to the central term in this sequence as the dual Atiyah bundle assoicated to $P_G$ relative to $T$.

In the situation where $X$ is a singular variety, we want a version of the dual Atiyah bundle which replaces Kähler differentials with Du Bois differentials. 
\begin{mdef}
 Let $0\to \Omega_{X/T}^{1,h}\to \mathrm{At}_{X/T}^{*,h}(P)\to \mathfrak{g}^*_P\to 0$ denote the sequence induced by the dual Atiyah sequence by pushing out along the map $\Omega^1_{X/T}\to \Omega^{1,h}_{X/T}$.
\end{mdef}

\subsubsection{Kähler forms on stacks} For a possibly singular stack $\mathcal{X}$, by a 1-form $\omega$ on $\mathcal{X}$ we mean a compatible system of Kähler $1$-forms $\omega_f\in H^0(T,\Omega^1_T)$ for all morphisms $f:T\to \mathcal{X}$, where $T$ is a smooth test object. If $\mathcal{X}$ is a smooth Artin stack, then $T^*\mathcal{X}$ admits a tautological $1$-form in this sense.

\section{The Tautological Du Bois Form on the Nilpotent Cone}

Let $\mathfrak{g}^*$ denote the dual Lie algebra of $G$ and $\mathcal{N}\subseteq \mathfrak{g}^*$ be the nilpotent cone. In this section we prove that the tautological 1-form $\omega_\tau$ on $T^*\mathcal{B}$ descends to a global Du Bois 1-form $\omega_\mathcal{N}$ on $\mathcal{N}$ along the Springer resolution
\[\sigma:T^*\mathcal{B}\to \mathcal{N}.\]

\begin{mprop}\label{descenting the tautological form}
    There exist a unique element $\omega_\mathcal{N}\in H^0(\mathcal{N},\Omega^{1,h}_\mathcal{N})$ such that $\sigma^*\omega_\mathcal{N}=\omega_\tau$ under the identification $H^0(T^*\mathcal{B},\Omega^{1,h}_{T^*\mathcal{B}})\cong H^0(T^*\mathcal{B},\Omega^1_{T^*\mathcal{B}})$
\end{mprop}

Consider the Steinberg variety
$\mathrm{St}_G=T^*\mathcal{B}\times_\mathcal{N}T^*\mathcal{B}$, which sits in a Cartesian square
\[
\begin{tikzcd}
\mathrm{St}_G \arrow[r, "\mathrm{pr}_1"] \arrow[d, "\mathrm{pr}_2"'] & T^*\mathcal{B} \arrow[d, "\sigma"] \\
T^*\mathcal{B} \arrow[r, "\sigma"']                                             & \mathcal{N}             
\end{tikzcd}\]
Since the Springer resolution is a proper surjective map, it suffices to prove that $\mathrm{pr}_1^*\omega_\tau=\mathrm{pr}_2^*\omega_\tau$ as Du Bois forms on the Steinberg variety.
\begin{mdef}
    Let $\omega_d=\mathrm{pr}_1^*\omega_\tau-\mathrm{pr}_2^*\omega_\tau\in H^0(\mathrm{St}_G,\Omega^1_{\mathrm{St}_G})$.
\end{mdef}

For each Weyl element $w\in W$, let $\mathrm{St}_G^w$ denote the locally closed subscheme of the Steinberg variety where the Borel subgroups have relative position $w$. Let $\overline{\mathrm{St}}_G^w$ denote the closure of $\mathrm{St}_G^w$. 

\begin{proof}[Proof of Proposition \ref{descenting the tautological form}]
By h-descent along the proper surjective map
\[\bigsqcup_{w\in W} \overline{\mathrm{St}}_G^w\to \mathrm{St}_G,\]
it suffices to prove that $\omega_d$ vanishes as a Du Bois form when pulled back to $\overline{\mathrm{St}}_G^w$ for all $w\in W$. Since the sheaf of Du Bois differential forms on any variety is torsion-free, it suffices to prove that $\omega_d$ vanishes as a Du Bois form when further pulled back to $\mathrm{St}_G^w$.  

The map $\mathrm{St}_G^w\to T^*\mathcal{B}\times T^*\mathcal{B}$ given by $(\mathrm{pr}_1,-\mathrm{pr}_2)$ identifies $\mathrm{St}_G^w$ with the conormal to the $G$-orbit $\mathcal{O}_w\subset \mathcal{B}\times \mathcal{B}$ consisting of pairs of Borel subgroups of relative position given by $w\in W$. Since $\mathrm{pr}_1^*\omega_\tau+\mathrm{pr}_2^*\omega_\tau$ is the tautological $1$-form on $T^*\mathcal{B}\times T^*\mathcal{B}\simeq T^*(\mathcal{B}\times \mathcal{B})$, it vanishes when restricted to the conormal of $\mathcal{O}_w$. Hence $\omega_d$ vanishes when pulled back to $\mathrm{St}_G^w$.
\end{proof}

To construct the map $H_0:\pi^{-1}\mathrm{Nilp}_{\mathrm{rel}}\to T^*\mathcal{M}_g$ from $\omega_\mathcal{N}$, we will need to understand how $\omega_\mathcal{N}$ transforms under the coadjoint action of $G$ on $\mathcal{N}$.

For a smooth variety $X$ with a $G$-action, let $\mathrm{act}_{T^*X}:G\times T^*X\to T^*X$ denote the map given by
\[G\times T^*X\overset{\mathrm{inv}_G\times\mathrm{id}}{\to} G\times T^*X\cong (G\times X)\times_XT^*X\to T^*(G\times X)\to T^*X.\]
where $(G\times X)\times_XT^*X\to T^*(G\times X)$ denotes the map given by differentiating the action map $G\times X\to X$. Let $m:G\times T^*X\to T^*G$ denote the composition
\[G\times T^*X\overset{\mathrm{inv}_G\times\mathrm{id}}{\to} G\times T^*X\cong (G\times X)\times_XT^*X\to T^*(G\times X)\to T^*G.\]

\begin{mlem}\label{actionsandmomentum}
Let $X$ be a smooth variety with a $G$-action and $G\times T^*X\to T^*X$ be the induced $G$-action on the cotangent space of $X$. Let $\omega_{T^*X}$ denote the tautological $1$-form on $T^*X$. Then
\[\mathrm{act}_{T^*X}^*\omega_{T^*X}=\mathrm{pr}_{T^*X}^*\omega_{T^*X}-m^*\omega_{T^*G},\]
where $m:G\times T^*X\to T^*G$ is the momentum map and $\omega_{T^*G}$ is the tautological $1$-form on $T^*G$. 
\end{mlem}

\begin{proof}
    The action map $G\times T^*X\to T^*X$ is given by the composition
    \[G\times T^*X\overset{\mathrm{inv}_G\times\mathrm{id}}{\to} G\times T^*X\cong (G\times X)\times_XT^*X\to T^*(G\times X)\to T^*X.\]
    Let $\omega_{T^*(G\times X)}$ be the tautological 1-form on $T^*(G\times X)$. Since the first map is given by differentiating the action map $G\times X\to X$, it follows that the pullback of $\omega_{T^*(G\times X)}$ to $G\times T^*X$ is equal to $\mathrm{pr}_{T^*X}^*\omega_{T^*X}$. Now $\omega_{T^*(G\times X)}$ is the sum of the pullback of $\omega_{T^*X}$ and $\omega_{T^*G}$ along the two projection maps. By definition the momentum map is given by the composition
    \[G\times T^*X\overset{\mathrm{inv}_G\times\mathrm{id}}{\to} G\times T^*X\cong (G\times X)\times_XT^*X\to T^*(G\times X)\to T^*G.\]
    Hence the statement follows.
\end{proof}

Identify $G\times \mathfrak{g}^*$ with $T^*G$ via left invariant 1-forms on $G$.

\begin{mprop}\label{p:transformation law}
    Let $\mathrm{act}_\mathcal{N}:G\times \mathcal{N}\to \mathcal{N}$ denote the coadjoint action. Then
    \[\mathrm{act}_\mathcal{N}^*\omega_\mathcal{N}=\mathrm{pr}_\mathcal{N}^* \omega_\mathcal{N}+ i^* \omega_G,\]
    where $i:G\times \mathcal{N}\to T^*G$ is the composition $G\times \mathcal{N}\subset G\times \mathfrak{g}^*\cong T^*G$ and $\omega_G$ is the tautological $1$-form on $T^*G$.
\end{mprop}
\begin{proof}
    Consider the commutative diagram
    \[
\begin{tikzcd}
G\times T^*\mathcal{B} \arrow[r, "\mathrm{act}_{T^*\mathcal{B}}"] \arrow[d, "{(\mathrm{id}_G,\sigma)}"'] & T^*\mathcal{B} \arrow[d, "\sigma"] \\
G\times \mathcal{N} \arrow[r, "\mathrm{act}_\mathcal{N}"]                                           & \mathcal{N}                  
\end{tikzcd}\]
Since 
\[m^*\omega_{T^*G}=-(\mathrm{id}_G,\sigma)^*i^*\omega_{T^*G},\]
it follows from Lemma \ref{actionsandmomentum} that 
\[(\mathrm{id}_G,\sigma)^*\mathrm{act}_\mathcal{N}^*\omega_\mathcal{N}=\mathrm{pr}_{T^*\mathcal{B}}^* \omega_\tau-m^*\omega_{T^*G}=(\mathrm{id}_G,\sigma)^*(\mathrm{pr}_\mathcal{N}^*\omega_\mathcal{N}+i^*\omega_{T^*G}).\]
Since $(\mathrm{id}_G,\sigma)$ is proper and surjective, we conclude that $\mathrm{act}_\mathcal{N}^*\omega_\mathcal{N}=\mathrm{pr}_\mathcal{N}^*\omega_\mathcal{N}+i^*\omega_{T^*G}$.
\end{proof}

We will also need the linearity of $\omega_\mathcal{N}$ under the scaling action of $\mathbb{G}_m$ on $\mathcal{N}$.

\begin{mprop}\label{p:linearity}
    Let $t$ be a coordinate on $\mathbb{G}_m$ and $w_1:\mathbb{G}_m\times \mathcal{N}\to \mathcal{N}$ denote the weight $1$-scaling action. Then
    \[w_1^* \omega_\mathcal{N}=t\cdot \mathrm{pr}_{\mathcal{N}}^* \omega_\mathcal{N}\]
    as Du Bois $1$-forms on $\mathbb{G}_m\times \mathcal{N}$.
\end{mprop}

\begin{proof}
    By arguing as in the previous proof and using the $\mathbb{G}_m$-equivariance of the Springer resolution, the statement follows from the linearity of the tautological $1$-form on $T^*\mathcal{B}$.
\end{proof}

\section{The Universal Nilpotent Cone}
In this section we explain how the tautological Du Bois form $\omega_\mathcal{N}$ on the nilpotent cone $\mathcal{N}$, constructed in the previous section, gives rise to a map $H_0:\pi^{-1}\mathrm{Nilp}_{\mathrm{rel}}\to T^*\mathcal{M}_g$, which in turn gives a construction and moduli interpretation of the universal nilpotent cone.

\subsubsection{Atiyah bundles and the cotangent space}
Recall that, for a $G$-bundle $P_G$ on a family of curves $X/T$, there is a canonical extension of vector bundles on $X$
\[0\to \Omega^1_{X/T}\to \mathrm{At}_{X/T}^*(P_G) \to \mathfrak{g}^*_{P_G}\to 0.\]
The absolute cotangent space $T^*\mathrm{Bun}_G(\mathcal{C}_g/\mathcal{M}_g)$ admits the following moduli interpretation:

A $T$-point on $T^*\mathrm{Bun}_G(\mathcal{C}_g/\mathcal{M}_g)$ is a triple $(X/T,P_G,v)$, where $X/T$ is a family of curves, $P_G$ is a $G$-bundle on $X$ and $v\in H^0(X,\mathrm{At}_{X/T}^*(P_G)\otimes \Omega^1_{X/T})$.

\noindent
We refer to the element $v$ as a \textit{absolute Higgs field} on $P_G$ relative to $T$. The projection $\pi$ from the absolute cotangent space to the cotangent space relative to $\mathcal{M}_g$ is given by sending an absolute Higgs field to its image under the map 
\[\mathrm{At}^*_{X/T}(P_G)\otimes \Omega^1_{X/T}\to \mathfrak{g}^*_{P_G}\otimes \Omega^1_{X/T}.\]
A relative Higgs field is equivalently described as a section of the map $\mathfrak{g}^*_{P_G}\times_{\mathbb{G}_m} \Omega^1_{X/T}\to X$ and is said to be relatively nilpotent if it factors through the closed subscheme
\[\mathcal{N}_{P_G}\times_{\mathbb{G}_m} \Omega^1_{X/T}\subset \mathfrak{g}^*_{P_G}\times_{\mathbb{G}_m} \Omega^1_{X/T}.\]
Or equivalently if, for all local trivializations of $P_G$ and $\Omega_{X/T}^1$ say on an étale neighborhood $U\to X$, the induced map $U\to \mathfrak{g}^*$ factors through $\mathcal{N}$.

\subsection{Constructing a section}
We will show that the Du Bois form $\omega_\mathcal{N}$ gives rise to a section $s$ of the map $\pi|_{\mathrm{Nilp}_\mathrm{rel}}: \pi^{-1}\mathrm{Nilp}_\mathrm{rel}\to \mathrm{Nilp}_\mathrm{rel}$.

\subsubsection{Description in local coordinates} Let us start by describing how to construct $s$ from $\omega_\mathcal{N}$ in local coordinates.

\begin{mrmk}
    The line bundle $\Omega^1_{X/T}$ occurs twice for different reasons. Once as the twist in Serre duality arising in the definition of Higgs bundles and once as the subbundle in the Atiyah sequence. In what follows, it will be useful to distinguish these two occurrences and replace the Serre twist with a general line bundle $L$.
\end{mrmk}

Let $X$ be a smooth variety of any dimension, $P_G$ a $G$-bundle on $X$ and $L$ a line bundle on $X$. A $L$-Higgs field $\overline{v}$ is a global section of $\mathfrak{g}^*_{P_G}\otimes L$. It is said to be nilpotent if for all local trivializations of $P_G$ and $L$, the induced map $\overline{v}:U\to \mathfrak{g}^*$ factors through $\mathcal{N}\subset \mathfrak{g}^*$. Given a nilpotent $L$-Higgs field $\overline{v}$, we want to produce a canonical global section $v$ of $\mathrm{At}^*_X(P_G)\otimes L$ lifting $\overline{v}$.

We start by producing local lifts given a choice of local trivializations of $P_G$ and $L$ and prove that these local lifts glue to a global section of $\mathrm{At}_X^*(P_G)\otimes L$. A local trivialization of $P_G$ over an étale neighborhood $U\to X$, provides a splitting of the dual Atiyah sequence
\[\Omega_X^1\to \Omega_X^1\oplus \mathfrak{g}^*\to \mathfrak{g}^*\]
as a vector bundle over $U$. Using the local trivialization of $L$, to specify a local lift of $\overline{v}$, we have to specify a local section of $\Omega_X^1$ over $U$. We then have to prove that these local differential forms satisfy the correct transformation property to define a global section of the twisted dual Atiyah bundle.

In the local trivialization of $P_G$ and $L$, the nilpotent $L$-Higgs field $\overline{v}$ gives rise to a map $\overline{v}:U\to \mathcal{N}$. Pulling back $\omega_\mathcal{N}$ along $\overline{v}$ and using the smoothness of $U$, we obtain a local section of $\Omega_X^1$ over $U$. Let $v$ be the local section of $\mathrm{At}^*_X(P_G)\otimes L$ given by $(\overline{v}^*\omega_\mathcal{N},\overline{v})\in H^0(U,\Omega_X^1\oplus \mathfrak{g}^*)$ under the chosen local trivializations.

To see that these local sections satisfy the correct transformation property, let $g:U\to G$ be a change of the local trivialization of $P_G$. The change of coordinates map on the dual Atiyah bundle is given by the matrix
\[\begin{pmatrix}
    \mathrm{id}_{\Omega_{X}^1}&m_g\\
    0& \mathrm{Ad}^*_g \\
\end{pmatrix},\]
where $m_g:\mathfrak{g}^*\to \Omega^1_U$ is given by considering $\mathfrak{g}^*$ as left-invariant differential forms on $G$ and pulling back along $g:U\to G$. The compatibility of the local section under this change of trivialization reduces to Proposition \ref{p:transformation law}. The compatibility under the change of local trivialization of $L$ reduces to the linearity of $\omega_\mathcal{N}$ under the $\mathbb{G}_m$-scaling action on $\mathcal{N}$.

To construct the desired section of $\pi^{-1}\mathrm{Nilp}_\mathrm{rel}\to \mathrm{Nilp}_\mathrm{rel}$, we assume $X/T$ is a family of curves, take $L=\Omega^1_{X/T}$ and send $v$ constructed above to its image under the map
\[\mathrm{At}^*_X(P_G)\otimes L\to \mathrm{At}^*_{X/T}(P_G)\otimes L.\]
Since this construction is functorial in maps between smooth test schemes, it defines the desired section
\[s:\mathrm{Nilp}_\mathrm{rel}\to \pi^{-1}\mathrm{Nilp}_\mathrm{rel}.\]

\subsubsection{Construction in the universal case} Alternatively to the local construction above, we can construct the section from the universal case. Let $\tau$ denote the diagonal global section of $\mathfrak{g}^*_ \mathcal{N}$ as a coherent sheaf on $\mathcal{N}$, $q:\mathcal{N}\to [\mathcal{N}/G]$ denote the quotient map, and $P_{G,\mathcal{N}}^u$ be the universal $G$-bundle on $[\mathcal{N}/G]$. Then $q^*P_{G,\mathcal{N}}^u$ comes with an isomorphism to the trivial $G$-bundle $P_G^0$ on $\mathcal{N}$. 

\begin{mlem}
The section $(\omega_\mathcal{N},\tau)\in H^0(\mathcal{N},\mathrm{At}_{\mathcal{N}}^{*,h}(P_G^0))$ descends to a global section $\eta_\mathcal{N}\in H^0([\mathcal{N}/G],\mathrm{At}_{[\mathcal{N}/G]}^{*,h}(P_{G,\mathcal{N}}^u))$.
\end{mlem}
\begin{proof}
Let $\mathrm{act}_\mathcal{N}:G\times \mathcal{N}\to \mathcal{N}$ denote the coadjoint action. To descent the section, it suffices to prove that the sections $\mathrm{act}_\mathcal{N}^* (\omega_\mathcal{N},\tau)$ and $\mathrm{pr}_\mathcal{N}^* (\omega_\mathcal{N},\tau)$ agree when considered as global sections of the Atiyah bundle of the trivial $G$-bundle on $G\times \mathcal{N}$. The descent datum for $P_{G,\mathcal{N}}^u$ induces the following map on dual Atiyah bundles over $G\times \mathcal{N}$
\[\begin{pmatrix}
    \mathrm{id}_{\Omega_{G\times\mathcal{N}}^{1,h}}& m\\
    0& \mathrm{Ad}^*_G \\
\end{pmatrix}\]
considered as an automorphism of $\mathrm{At}_{\mathcal{N}\times G}^{*,h}(P_G^0)\simeq \Omega_{\mathcal{N}\times G}^{1,h}\oplus \mathfrak{g}^*_{G\times \mathcal{N}}$ as a coherent sheaf on $G\times \mathcal{N}$. Here $m:\mathfrak{g}^*_{ G\times \mathcal{N}}\to \Omega_{G\times \mathcal{N}}^{1,h}$ is given by considering $\mathfrak{g}^*$ as left invariant 1-forms on $G$ and pulling them back along the  projection $G\times \mathcal{N}\to G$ to obtain global sections of $\Omega^{1,h}_{G\times \mathcal{N}}$. 

Since $\tau$ is the diagonal section of $\mathfrak{g}^*_\mathcal{N}$, it follows that $m$ sends the pullback of $\tau$ to the $1$-form on $G\times \mathcal{N}$ given by pulling back $\omega_{T^*G}$ along the map $G\times \mathcal{N}\to \mathfrak{g}^*\times G\simeq T^*G$. Thus the compatibility of $(\omega_\mathcal{N},\tau)$ with the descent datum reduces to Proposition \ref{p:transformation law}.
\end{proof}

Let $X$ be a smooth variety, $L$ a line bundle on $X$ and $P_G$ a $G$-bundle on $X$. Using the coadjoint action of $G$ on $\mathcal{N}$, we can form the twisted product $\mathcal{N}_{P_G}$. Since coadjoint action of $G$ on $\mathcal{N}$ commutes with the weight $1$ scaling action of $\mathbb{G}_m$ on $\mathcal{N}$, we can further twist this by the $\mathbb{G}_m$-bundle defined by $L$ to obtain $Y:=\mathcal{N}_{P_G}\times_{\mathbb{G}_m} L$.

Let $P_{G,Y}^u$ be the universal $G$-bundle on $\mathcal{N}_{P_G}\times_{\mathbb{G}_m} L$ and $\overline{v}_{P_G}\in H^0(\mathcal{N}_{P_G}\times_{\mathbb{G}_m} L,\mathfrak{g}^*_{P_{G,Y}^u}\otimes L)$ be the canonical global section.

\begin{mlem}
    The global section $\eta_\mathcal{N}\in H^0([\mathcal{N}/G],\mathrm{At}_{[\mathcal{N}/G]}^{*,h}(P_{G,\mathcal{N}}^u))$ scales with weight $1$ under the $\mathbb{G}_m$-action on $[\mathcal{N}/G]$ and hence gives rise to a global section 
    \[v_{P_G,L}\in H^0(\mathcal{N}_{P_G}\times_{\mathbb{G}_m} L,\mathrm{At}_{Y}^{*,h}(P_{G,Y}^u)\otimes L)\]
    lifting $\overline{v}_{P_G}$.
\end{mlem}

\begin{proof}
    This follows from Proposition \ref{p:linearity} and the linear scaling of the diagonal section $\tau$ of $\mathfrak{g}^*_\mathcal{N}$.
\end{proof}

For a family of curves $X\to T$ over a smooth test scheme $T$, let $L=\Omega^1_{X/T}$. Given a $G$-bundle $P_G$, we denote the resulting section $v_{P_G,L}$ simply by $v_{P_G}$.

Given a $T$-valued point $(X/T,P_G,\overline{v})$ of $\mathrm{Nilp}_\mathrm{rel}$, we obtain a section $s_{\overline{v}}:X\to \mathcal{N}_{P_G}\times_{\mathbb{G}_m} \Omega^1_{X/T}$ such that $\overline{v}=s_{\overline{v}}^*(\overline{v}_{P_G})$. Thus pulling back $v_{P_G}$ along $s_{\overline{v}}$ gives a global section of $ H^0(X,\mathrm{At}^*_X(P_G)\otimes \Omega^1_{X/T})$, let $v\in H^0(X,\mathrm{At}^*_{X/T}(P_G)\otimes \Omega^1_{X/T})$ denote its image in the twisted dual Atiyah bundle relative to $T$. Then $v$ provides a  lift of $\overline{v}$. Thus sending $(X/T,P_G,\overline{v})$ to the $T$-point $(X/T,P_G,v)$ of $\pi^{-1}\mathrm{Nilp}_\mathrm{rel}$, we have produced the section 
\[s:\mathrm{Nilp}_\mathrm{rel}\to \pi^{-1}\mathrm{Nilp}_\mathrm{rel}.\]

\subsubsection{The universal nilpotent cone} Consider the map
\[\pi^{-1}\mathrm{Nilp}_\mathrm{rel}\overset{\mathrm{id}-s\circ\pi}{\to} T^*\mathrm{Bun}_G(\mathcal{C}_g/\mathcal{M}_g).\]
Since $s$ is a section of $\pi|_{\mathrm{Nilp}_\mathrm{rel}}$, this factors through $\mathrm{Bun}_G(\mathcal{C}_g/\mathcal{M}_g)\times_{\mathcal{M}_g}T^*\mathcal{M}_g$ to give a map
\[H_0:\pi^{-1}\mathrm{Nilp}_\mathrm{rel}\to T^*\mathcal{M}_g.\]
\begin{mdef}
    Let $\mathrm{Nilp}_u$ denote the closed substack given by $H_0^{-1}(\mathbf{0}_{T^*\mathcal{M}_g })$, where $\mathbf{0}_{T^*\mathcal{M}_g}$ denotes the closed substack of $T^*\mathcal{M}_g$ given by the zero section.
\end{mdef}

\begin{mrmk}
    We say that a map of stacks $\mathcal{Z}\to \mathcal{X}$ is a closed immersion, if for all smooth test objects $T\to \mathcal{X}$ the pullback $\mathcal{Z}\times_\mathcal{X} T\to T$ is representable as a stack by a closed subscheme of $T$.
\end{mrmk}

To relate this construction of the universal nilpotent cone to $B$-reductions, define $\pi^{-1}\mathrm{Nilp}_{\mathrm{rel},B}$ to be the pullback in the following diagram
\[
\begin{tikzcd}
{\pi^{-1}\mathrm{Nilp}_{\mathrm{rel},B}} \arrow[d] \arrow[r]   & T^*\mathrm{Bun}_G\times_{\mathrm{Bun}_G}\mathrm{Bun}_B \arrow[d] \\
T^*\mathcal{M}_g\times_{\mathcal{M}_g}\mathrm{Bun}_B \arrow[r] & T^*\mathrm{Bun}_B                                               
\end{tikzcd}\]
The composition $\pi^{-1}\mathrm{Nilp}_{\mathrm{rel},B}\to T^*\mathrm{Bun}_G\times_{\mathrm{Bun}_G}\mathrm{Bun}_B\to T^*\mathrm{Bun}_G$ factors through a map $\mathrm{pr}:\pi^{-1}\mathrm{Nilp}_{\mathrm{rel},B}\to\pi^{-1}\mathrm{Nilp}_\mathrm{rel}$. Let $H_{0,B}:\pi^{-1}\mathrm{Nilp}_{\mathrm{rel},B}\to T^*\mathcal{M}_g$ denote the composition of the left vertical morphism in the diagram above followed by the projection to $T^*\mathcal{M}_g$. 

\begin{mprop}
    The map $H_0$ is uniquely determined by the property $H_0\circ \mathrm{pr}\simeq H_{0,B}$.
\end{mprop}

\begin{proof}
    Since we are only considering smooth test objects and $\mathrm{pr}:\pi^{-1}\mathrm{Nilp}_{\mathrm{rel},B}\to\pi^{-1}\mathrm{Nilp}_\mathrm{rel}$ is surjective on geometric points, the uniqueness follows once we establish the equivalence. 

    A $T$-valued point of $\pi^{-1}\mathrm{Nilp}_{\mathrm{rel},B}$ is given by a family of curves $X/T$, a $B$-bundle $P_B$, an absolute Higgs field $v$ on the induced $G$-bundle $P_G$, such that $v$ is in the kernel of the composition
    \[H^0(X,\mathrm{At}^*_{X/T}(P_G)\otimes \Omega^1_{X/T})\to H^0(X,\mathrm{At}^*_{X/T}(P_B)\otimes \Omega^1_{X/T})\to H^0(X,\mathfrak{b}^*_{P_B}\otimes \Omega^1_{X/T}).\]
    In other words, the image of $v$ in $H^0(X,\mathrm{At}^*_{X/T}(P_B)\otimes \Omega^1_{X/T})$ is the image of a quadratic differential $w\in H^0(X,\Omega^{\otimes 2}_{X/T})$ under the map $H^0(X,\Omega^{\otimes 2}_{X/T})\to H^0(X,\mathrm{At}^*_{X/T}(P_B)\otimes \Omega^1_{X/T})$. By construction the quadratic differential $w$ considered as a $T$-point of $T^*\mathcal{M}_g$ is the image under the map $H_{0,B}$. Here we identify $T^*\mathcal{M}$ is the moduli space of a curve with a quadratic differential relative to the base.

    On the other hand, étale locally on $U\to X$ choose trivializations of $P_G$ and $L=\Omega^1_{X/T}$. Given such trivializations, $v$ and the $B$-reduction $P_B$ define a map $U\to T^*\mathcal{B}$. Pulling back the tautological $1$-form on $T^*\mathcal{B}$ gives a $1$-form on $U$, which we can identify with a quadratic form under the chosen trivialization of $L$. The quadratic differential $w$ restricted to $U$ is exactly the quadratic component of $v$ restricted to $U$ with respect to the chosen trivializations minus the quadratic form obtained from pulling back the tautological 1-form on $T^*\mathcal{B}$. Since $\omega_\mathcal{N}$ descends the tautological $1$-form along the Springer resolution, it follows that the pullback of the tautological $1$-form on $T^*\mathcal{B}$ is equal to the pullback of $\omega_\mathcal{N}$ under the map $U\to \mathcal{N}$. Thus the restriction of $w$ to $U$ is equal to the quadratic component of $v$ under the chosen trivialization minus the quadratic form obtained from the pullback of $\omega_\mathcal{N}$, which is exactly the restriction of the quadratic differential defined by $H_0$ to $U$. By covering $X$ with such $U$, we conclude that the two maps agree.
\end{proof}

We now establish the fundamental properties of $\mathrm{Nilp}_u$ using $H_0$. Since $H_0$ is equivariant for  the $\mathbb{G}_m$-action, it follows that $\mathrm{Nilp}_u$ is a cone. Note that the map $\pi$ restricts to an equivalence $\mathrm{Nilp}_u\to \mathrm{Nilp}_\mathrm{rel}$ with inverse given by $s$. Thus the universal nilpotent cone and the relative nilpotent cone are equivalent as stacks on smooth test varieties. The last remaining property we will prove is that $\mathrm{Nilp}_u$ is Lagrangian, this will follow from the next proposition. Let $\omega_{T^*\mathcal{M}_g}$ denote the tautological $1$-form on $T^*\mathcal{M}_g$.

\begin{mprop}\label{H0and1forms}
    The $1$-form on $\pi^{-1}\mathrm{Nilp}_\mathrm{rel}$ given by the pullback of $\omega_{T^*\mathcal{M}_g}$ along $H_0$ agrees with the restriction of the tautological $1$-form $\omega_{T^*\mathrm{Bun}_G}$ on $T^*\mathrm{Bun}_G(\mathcal{C}_g/\mathcal{M}_g)$.
\end{mprop}

\begin{proof}
    We start by proving that the $1$-forms agree after pullback along the map $\pi^{-1}\mathrm{Nilp}_{\mathrm{rel},B}\to \pi^{-1}\mathrm{Nilp}_{\mathrm{rel}}$. Consider the following diagram
\[
\begin{tikzcd}
\pi^{-1}\mathrm{Nilp}_{\mathrm{rel}} \arrow[r]                            & T^*\mathrm{Bun}_G                                                          \\
{\pi^{-1}\mathrm{Nilp}_{\mathrm{rel},B}} \arrow[r] \arrow[u] \arrow[d]    & T^*\mathrm{Bun}_G\times_{\mathrm{Bun}_G}\mathrm{Bun}_B \arrow[u] \arrow[d] \\
T^*\mathcal{M}_g\times_{\mathcal{M}_g} \mathrm{Bun}_B \arrow[r] \arrow[d] & T^*\mathrm{Bun}_B                                                          \\
T^*\mathcal{M}_g                                                          &                                                                           
\end{tikzcd}\]
The composition of $\pi^{-1}\mathrm{Nilp}_{\mathrm{rel},B}\to \pi^{-1}\mathrm{Nilp}_{\mathrm{rel}}$ with $H_0$ is $H_{0,B}$, which by definition is the composition of the two left-most downward morphisms. The pullback of $\omega_{T^*\mathcal{M}_g}$ to $T^*\mathcal{M}_g\times_{\mathcal{M}_g} \mathrm{Bun}_B$ is equal to the pullback of the tautological $1$-form $\omega_{T^*\mathrm{Bun}_B}$ under the map to $T^*\mathrm{Bun}_B$. The pullback of $\omega_{T^*\mathrm{Bun}_B}$ to $T^*\mathrm{Bun}_G\times_{\mathrm{Bun}_G}\mathrm{Bun}_B$ is equal to the pullback of $\omega_{T^*\mathrm{Bun}_G}$ under the first projection. Thus the pullback of $\omega_{T^*\mathcal{M}_g}$ to $\pi^{-1}\mathrm{Nilp}_{\mathrm{rel},B}$ is equal to the pullback of $\omega_{T^*\mathrm{Bun}_G}$.

The map $\pi^{-1}\mathrm{Nilp}_{\mathrm{rel},B}\to \pi^{-1}\mathrm{Nilp}_{\mathrm{rel}}$ is surjective on geometric points. For any smooth test scheme $T\to \pi^{-1}\mathrm{Nilp}_\mathrm{rel}$, we have to prove that the 1-forms agree after pullback to $T$. The map $T\times_{\pi^{-1}\mathrm{Nilp}_\mathrm{rel}} \pi^{-1}\mathrm{Nilp}_{\mathrm{rel},B}\to T$ is surjective on geometric points and the two 1-forms agree after further pullback along this map. It follows by generic smoothness that the two forms are equal on a dense open subset of $T$. Since $T$ is smooth we conclude that the two $1$-forms are equal when pulled back to $T$. 
\end{proof}

\begin{mcor}
The universal nilpotent cone $\mathrm{Nilp}_u$ is Lagrangian in $T^*\mathrm{Bun}_G$.
\end{mcor}

\begin{proof}
By Proposition \ref{H0and1forms} the tautological 1-form on $T^*\mathrm{Bun}_G$ restricted to $\mathrm{Nilp}_u$ equal to the pullback along $H_0$ of the tautological $1$-form on $T^*\mathcal{M}_g$. Since the tautological $1$-form on $T^*\mathcal{M}_g$ vanishes on the zero-section $0_{T^*\mathcal{M}_g}$, it vanishes when pulled back to $\mathrm{Nilp}_u$. Thus $\mathrm{Nilp}_u$ is isotropic.

Knowing that $\mathrm{Nilp}_u$ is isotropic, all that remains is to prove that the irreducible components of $\mathrm{Nilp}_u$ have dimension at least
\[\operatorname{dim} \mathrm{Bun}_G(\mathcal{C}_g/\mathcal{M}_g).\]
We give a proof assuming $g>1$, for the general case see \cite[Lemma 4.2.1]{nadler2026global} . Since $\mathrm{Nilp}_u$ maps isomorphically to $\mathrm{Nilp}_\mathrm{rel}$, we are reduced to proving the analogous statement for the relative nilpotent cone. By definition, the relative nilpotent cone is the preimage of the zero-section under the Hitchin fibration
\[T^*_{\mathcal{M}_g}\mathrm{Bun}_G(\mathcal{C}_g/\mathcal{M}_g)\to \mathrm{Hitch}_G.\]
Let $C$ be the dimension of the center of $G$. The irreducible components of the stack $T^*_{\mathcal{M}_g}\mathrm{Bun}_G(\mathcal{C}_g/\mathcal{M}_g)$ have dimension $\ge 2\operatorname{dim} \mathrm{Bun}_G+C-\operatorname{dim} \mathcal{M}_g$ and $\mathrm{Hitch}_G$ is smooth over $\mathcal{M}_g$ and hence, by the analogous statement for a single curve, has dimension equal to $\operatorname{dim} \mathrm{Bun}_G+C$.  Thus the irreducible components of the relative nilpotent cone has dimension at least 
\[(2\operatorname{dim} \mathrm{Bun}_G+C-\operatorname{dim} \mathcal{M}_g)-(\operatorname{dim} \mathrm{Bun}_G+C)+\operatorname{dim} \mathcal{M}_g=\operatorname{dim} \mathrm{Bun}_G\]
\end{proof}

\subsubsection{The moduli description} As promised this construction of the $\mathrm{Nilp}_u$ gives a clear moduli interpretation of $T$-valued points for smooth test schemes $T$. To make this completely explicit, we state a very concrete version of the moduli description.

For $T$ a smooth variety, a $T$-valued point of $\mathrm{Nilp}_u$ is given by $(X/T,P_G,v)$ such that $v$ is relatively nilpotent and for any local trivialization of $P_G$ and $\Omega^1_{X/T}$ over an étale neighborhood $U\to X$, the quadratic differential component of $v$ under the trivializations is equal to the pullback of $\omega_\mathcal{N}$ along the induced map $U\to \mathcal{N}$, considered as a quadratic differential form on $X$ relative to $T$ using the given trivialization of $\Omega_{X/T}^1$.

\begin{mrmk}
    The reader might object that the above characterization is not intrinsic as it refers to trivializations of $P_G$ and $\Omega^1_{X/T}$. However, note that the above is a condition on any trivialization, whereas the characterization in terms of $B$-reductions is a condition about the existence of a very special $B$-reduction. The transformation properties ensure that the condition does not depend on the local trivializations.
\end{mrmk}

\section{Extension to $\operatorname{Spec }\mathbb{Z}$ for $\mathrm{PGL}_2$}

We end the paper by an example which indicate that our construction of the universal nilpotent cone is robust enough to work in far greater generality than just base fields of characteristic $0$. Working over $\operatorname{Spec } \mathbb{Z}$ with $G=\mathrm{PGL}_{2,\mathbb{Z}}$, we prove that the tautological $1$-form on $T^*\mathcal{B}$ in some sense descends to a global $1$-form on $\mathcal{N}$. More precisely, we prove the following proposition.

Let $\omega_\mathrm{reg}$ denote the 1-form on $\mathcal{N}_\mathrm{reg}$ obtained from the tautological $1$-form on $T^*\mathcal{B}$ and the isomorphism given by the Springer resolution on the regular locus of $\mathcal{N}$.

\begin{mprop}\label{WorksoverZ}
To every pair of a locally factorial scheme $X$ with a map $N:X\to \mathcal{N}$, there is a unique assignment of a global Kähler $1$-form $\omega_N$ on $X$, which is functorial in $(X,N)$ and extends the pullback of $\omega_\mathrm{reg}$ on the preimage of regular locus $\mathcal{N}_\mathrm{reg}$ under $N$.
\end{mprop}

From this proposition we could follow the same construction as before to give a definition of $\mathrm{Nilp}_u$ defined over $\operatorname{Spec }\mathbb{Z}$, however, we will not write this construction out in detail.

The construction in Proposition \ref{WorksoverZ} is carried out in terms of \textit{factorization systems}, rather than simply checking that for $N$ not identically $0$ the $1$-form $N^*\omega_\mathrm{reg}$ extends to all of $X$, to ensure that the extensions satisfy functoriality in the pair $(X,N)$.

To begin the construction, let $N$ be given by a matrix
\[\begin{pmatrix}
    x& y\\z&-x
\end{pmatrix}\]
where $x,y,z$ are global functions on $X$ satisfying $x^2+yz=0$. On the regular locus of $\mathcal{N}$, the $1$-form $\omega_{\mathcal{N}_{reg}}$ is given in coordinates by the formula $-dx+x\operatorname{dlog} y=dx-x\operatorname{dlog} z$.

We want to make sense of the global $1$-form $x\operatorname{dlog}y$ on $X$, to do this in a functorial way we introduce the notion of a factorization system. Locally we may assume that $X=\operatorname{Spec }R$ for an integral domain $R$. To motivate the next definition, assume momentarily that $R$ is a unique factorization domain $R$. The equation $x^2+yz=0$ implies that $y|x^2$. Since $R$ is a UFD, if $y$ is irreducible this implies that $y|x$. In particular we can choose a factorization $y=y_1\cdots y_r$ such that $x=g_iy_i$. This motivates the following definition.

\begin{mdef}
     Let $R$ be an integral domain and $x,y \in R$ be elements such that $y|x^2$. A \textit{factorization system} for $y$ in $x$ is a factorization $y=y_1\cdots y_r$, for $y_i\in R$, together with $g_i\in R$ such that $x=g_iy_i$.
\end{mdef}

Note that we do not assume that $y_i$ is irreducible. Given such a factorization system for $y$ in $x$, we define a global $1$-form on $X$ by
\[x\operatorname{dlog}y:=\sum_{i=1}^r g_idy_i.\]
The following lemma justifies that the symbol $x\operatorname{dlog}y$ does not refer to the chosen factorization system for $y$ in $x$.
\begin{mlem}
For $R$ locally factorial, the $1$-form $\sum_{i=1}^r g_idy_i$ is independent of the chosen factorization system
\end{mlem}

\begin{proof}
Since the statement can be checked on each stalk, we may assume that $R$ is a UFD. If $x=0$, then for each $i$ either $g_i=0$ or $y_i=0$ and hence the resulting 1-form is $0$ and in particular independent of the factorization system. If $y=0$, then $x=0$, so we may assume $y\neq 0$. By factoring $y$, it follows by the Leibniz rule and unique factorization, that the $1$-form associated to any factorization system is equal to the 1-form determined by the factorization system corresponding to the unique factorization of $y$.
\end{proof}

Since $X$ is assumed to be a locally factorial scheme, the stalks are unique factorization domains. Given any $x,y\in R$ for a unique factorization domain $R$ such that $y|x^2$, we can always find a factorization system for $y$ in $x$, by factoring $y$ into irreducibles. For a locally factorial integral scheme $X$ with a point $p\in X$ and global sections $x,y$ such that $y|x^2$, we can find a  factorization system $y_1,\dots y_r$ for $y$ in $x$ inside the stalk of $X$ at $p$. By spreading out, we find a factorization system for $y$ in $x$ defined in a Zariski neighborhood $U$ of $p$. Thus we can define the $1$-form $x\operatorname{dlog} y$ as a $1$-form on $U$. By the independence of the factorization system, we see that as $U$ ranges over a Zariski cover of $X$ where a factorization system can be chosen on each open set, the resulting 1-forms $x\operatorname{dlog} y$ are compatible and descent to a global $1$-form on $X$, which we denote by the same symbol.

Thus given $(X,N)$ we have constructed the desired $1$-form $\omega_N=-dx+x\operatorname{dlog} y$, which extends $N^*\omega_\mathrm{reg}$. All that remains is to prove functoriality. Having set up things in terms of factorization systems, functoriality follows immediately from the trivial lemma below.

\begin{mlem}
     Let $f:R\to R'$ be a morphism of integral domains, $x,y\in R$ such that $y|x^2$, and $y_1,\dots, y_r$ be a factorization system for $y$ in $x$. Then $f(y_1),\dots, f(y_r)$ is a factorization system for $f(y)$ in $f(x)$. 
\end{mlem}
 
\subsubsection{Uniqueness} Let $\mathbb{A}^3_\mathbb{Z}=\operatorname{Spec} \mathbb{Z}[r,s,t]$. Consider the map $N:\mathbb{A}^3_\mathbb{Z}\to \mathcal{N}$ given by the nilpotent matrix
\[
\begin{pmatrix}
    rst& rs^2\\
    -rt^2 & -rst
\end{pmatrix}
\]
Since the preimage of the regular locus in $\mathcal{N}$ is dense in $\mathbb{A}^3_\mathbb{Z}$ and the module of differential forms on $\mathbb{A}^3_\mathbb{Z}$ is free, the differential form associated to $(\mathbb{A}^3_\mathbb{Z},N)$ is uniquely determined by the compatibility with $\omega_\mathrm{reg}$ on the regular locus. Let $R$ be a unique factorization domain, then any map $\operatorname{Spec} R\to \mathcal{N}$ factors through $N:\mathbb{A}^3_\mathbb{Z}\to \mathcal{N}$. To see this claim, if $y=0$ take $r=-z$, $s=0$ and $t=1$. If $y\neq 0$, let $d=\operatorname{gcd}(x,y)$ and write $x=td$ and $y=sd$. The equation implies that $s|d$, let $d=rs$. This gives the desired factorization. Hence the form on $\operatorname{Spec}R$ is uniquely determined by the functoriality condition. 

For a general locally factorial scheme, it suffices to prove that all stalks are uniquely determined. By the functoriality condition this reduces to the case of the spectrum of a unique factorization domain.

\printbibliography

@article{nadler2021automorphic,
  title={Automorphic gluing functor in Betti geometric Langlands},
  author={Nadler, David and Yun, Zhiwei},
  journal={arXiv preprint arXiv:2105.12318},
  year={2021}
}

@article{nadler2026global,
  title={{The global nilpotent cone for universal curves}},
  author={Nadler, David and Yun, Zhiwei},
  journal={arXiv preprint arXiv:2603.01261},
  year={2026}
}

@misc{huber2014differentialformshtopology,
      title={Differential forms in the h-topology}, 
      author={Annette Huber and Clemens Jörder},
      year={2014},
      eprint={1305.7361},
      archivePrefix={arXiv},
      primaryClass={math.AG},
      url={https://arxiv.org/abs/1305.7361}, 
}
\end{document}